\documentclass[journal]{IEEEtran}

\usepackage{amssymb}
\usepackage[ruled,vlined,linesnumbered]{algorithm2e}
\usepackage{amsthm}
\usepackage{amsfonts}
\usepackage{amsmath}
\usepackage{booktabs}
\usepackage{bm}
\usepackage{color}
\usepackage{cite}
\usepackage{comment}
\usepackage{diagbox}
\usepackage{epsfig}
\usepackage{euscript}
\usepackage{fancybox}
\usepackage{graphicx}
\usepackage{hyperref}
\usepackage{lipsum}
\usepackage{pgf}
\usepackage{psfrag}
\usepackage{pifont}
\usepackage{mathrsfs}
\usepackage{multirow}
\usepackage{setspace}
\usepackage{stfloats}
\usepackage{tabularx}
\usepackage{textcomp}

\newcolumntype{C}[1]{>{\centering\arraybackslash}p{#1}}

\newtheorem{Theorem}{Theorem}
\newtheorem{Corollary}{Corollary}

\newtheorem{Definition}{Definition}

\SetKwFunction{Dothat}{Do that}

    \def\Complex{{\rm\rule[.23ex]{.03em}{1.1ex}\kern-.3em{C}}}

    \newcommand{\be}{\begin{equation}} \newcommand{\ee}{\end{equation}}
    \newcommand{\bea}{\begin{eqnarray}} \newcommand{\eea}{\end{eqnarray}}
    \newcommand{\benum}{\begin{enumerate}} \newcommand{\eenum}{\end{enumerate}}

    \newcommand*{\argmax}{\operatornamewithlimits{argmax}\limits}

\begin{document}

\title{Distributed Deep Reinforcement Learning for Intelligent Load Scheduling in Residential Smart Grids}

\author{ \IEEEauthorblockN{Hwei-Ming Chung, \IEEEmembership{Student Member, IEEE}, Sabita Maharjan, \IEEEmembership{Senior Member, IEEE}, Yan Zhang, \IEEEmembership{Fellow, IEEE}, and Frank Eliassen, \IEEEmembership{Member, IEEE} }

\thanks{This work was supported by Norwegian Research Council under Grants 275106 (LUCS project), 287412 (PACE project), and 267967 (SmartNEM project).
}

\thanks{H.-M. Chung and F. Eliassen are with the Department of Informatics, University of Oslo, Oslo 0373, Norway, (e-mail: hweiminc@ifi.uio.no, frank@ifi.uio.no).}

\thanks{S. Maharjan and Y. Zhang are with Department of Informatics, University of Oslo, Oslo 0373, Norway; and Simula Metropolitan Center for Digital Engineering, Oslo 0167, Norway. (e-mail: yanzhang@ieee.org, sabita@ifi.uio.no).
}

}

\maketitle

\begin{abstract}
The power consumption of households has been constantly growing over the years.
To cope with this growth, intelligent management of the consumption profile of the households is necessary, such that the households can save the electricity bills, and the stress to the power grid during peak hours can be reduced. 
However, implementing such a method is challenging due to the existence of randomness in the electricity price and the consumption of the appliances.
To address this challenge, we employ a model-free method for the households which works with limited information about the uncertain factors.
More specifically, the interactions between households and the power grid can be modeled as a non-cooperative stochastic game, where the electricity price is viewed as a stochastic variable.
To search for the Nash equilibrium (NE) of the game, we adopt a method based on distributed deep reinforcement learning.
Also, the proposed method can preserve the privacy of the households.
We then utilize real-world data from Pecan Street Inc., which contains the power consumption profile of more than $1,000$ households, to evaluate the performance of the proposed method.
In average, the results reveal that we can achieve around $12\%$ reduction on peak-to-average ratio (PAR) and $11\%$ reduction on load variance.
With this approach, the operation cost of the power grid and the electricity cost of the households can be reduced.

\end{abstract}

\section{Introduction}
\IEEEPARstart{W}{ith} a growing population, an increasing standard of living, and more power-demanding appliances being used in households, the energy consumption in residential households has increased considerably in the recent years, and is anticipated to grow even further \cite{2019-EIA-report}.
Improvements in energy efficiency have not been significant enough to counteract the increasing demand \cite{2018-EEA-report}.
To overcome this situation, the deployment of intelligent devices and communication infrastructure in smart grids becomes an important initiative.
By doing so, the demand side is able to play an active role in energy management to balance demand and supply.
More specifically, the demand side can change the consumption profile based on information (such as electricity price and generation capacity) from the supply side.
For example, the authors in \cite{2011-DR-first} proposed several concepts for scheduling the consumption of household appliances to reduce the electricity cost.

Subsequently, many research papers \cite{2013-DR-sabita,2016-DR-sabita,2013-DR-ICB,2016-DR-industrial-load,2018-Eh-DR-incentive,2017-DR-privacy,2018-DR-uncertainty,2018-pil-DR-QOE} have suggested several mechanisms to schedule the consumption of the appliances.
The authors in \cite{2013-DR-sabita,2016-DR-sabita} applied game theory to model the interaction between the utility companies and the customers to reduce the power consumption.
A real-time pricing scheme was adopted in \cite{2013-DR-ICB}, and a genetic algorithm was utilized to minimize the electricity cost.
Instead of scheduling the appliances in the residential area, the scheduling of industrial loads was considered in \cite{2016-DR-industrial-load}.
Then, an incentive scheme was applied in \cite{2018-Eh-DR-incentive} to encourage more households to participate in load scheduling.
When scheduling appliances, privacy-sensitive information of households may be revealed to the utility companies.
In this regard, the solutions proposed in \cite{2013-DR-sabita,2016-DR-sabita,2013-DR-ICB,2016-DR-industrial-load,2018-Eh-DR-incentive} may be of limited applicability.
Therefore, the authors in \cite{2017-DR-privacy} proposed an algorithm adopting a randomized approach to address the privacy issue while scheduling household appliances.
On the other hand, the power system will include even more renewable energy resources in the future.
Therefore, it is important to incorporate the uncertainty associated with renewable energy resources for scheduling the consumption of appliances as stated in \cite{2018-DR-uncertainty}.
The authors in \cite{2018-pil-DR-QOE} further considered the quality of experience (QoE) when scheduling the load with renewable energy.

However, with more and more households connected to the power grid, the computational complexity of the centralized algorithms has become a significant issue.
Therefore, several works have proposed distributed algorithms for scheduling the appliances \cite{2014-DR-dis,2018-DR-dis, 2018-DR-dis-shape, 2018-DR-dis-EV,2019-genn-DR-emission}.
In \cite{2014-DR-dis}, the distributed algorithm required the coordination between households. 
The authors in \cite{2018-DR-dis} proposed a scheme where the households can have access to the energy market, and the independent system operator (ISO) can send the incentive signal to the households such that the generation cost can be minimized.
The authors in \cite{2018-DR-dis-shape} introduced a primal-dual method to design a distributed algorithm, {so that the necessary amount of communication between households can be reduced. 
In \cite{2018-DR-dis-EV}, a distributed algorithm was designed based on the architecture in \cite{2018-DR-dis} where households can directly participate in the energy market.
The distributed algorithm jointly scheduled the charging behavior of electric vehicles (EVs) to reduce the procurement cost from the energy market.
Then, distributed load scheduling was applied in commercial ports in \cite{2019-genn-DR-emission}.

In the above works, the scheduling problem is mainly solved with model-based control approaches.
These approaches provided good results for scheduling the consumption of the appliances.
However, detailed information on the system model is required as input to the model-based control approaches.
For example, the transition probability of the electricity price should be modeled; however, electricity price is generated real-time, and therefore it cannot be obtained in advance.
Obtaining such an accurate model and maintaining its accuracy are non-trivial tasks.
In this context, model-free control approaches are considered as valuable alternatives to model-based control solutions.
An advantage of model-free methods is that they can obtain similar performance as model-based methods without relying on the knowledge of the system.
The most popular model-free control paradigm is reinforcement learning (RL).
RL has also been deployed in several works in the energy domain.
For instance, the authors in \cite{2009-ML-MDP-compare} applied both a model-based method and RL to control the power oscillation damping.
The results reveal that RL can be more robust compared to the model-based method.  
Then, the problem of scheduling the consumption of the appliances can be formulated with a Markov decision process (MDP) and solved by applying traditional RL techniques (e.g., in \cite{2015-RL-DR} and \cite{2019-lu-DR-RL}, the authors deployed Q-learning).
The authors in \cite{2018-water-heater-RL-DR} focused on using Q-learning to schedule the power consumption of water heaters.
A similar framework was discussed in \cite{2016-RL-DR-kim}, where the authors adopted distributed Q-learning to learn the scheduling policy.
The authors in \cite{2017-RL-DR} focused on scheduling the load profile of thermostatically controlled loads (e.g., water heaters or heat pumps).
Load scheduling with RL was incorporated with privacy-preserving mechanism in \cite{2019-amir-DR-RL-privacy}.
However, Q-learning needs to discretize the state and the action spaces to build the Q-table; the performance relies heavily on the resolution of the discretization.
Also, higher resolution of the discretization increases the computation time.

Recently, compared with the model-based methods, model-free methods have obtained remarkable breakthroughs.
The researchers from OpenAI found out that deep neural network (DNN) can be a good approximator to generate the Q-table; thus, the discretization of the state space in the traditional RL techniques is not necessary.
This framework is called Deep Q-Network (DQN) \cite{2015-DQN-nature}, and it has been applied in \cite{2018-building-RL,2018-EV-charg-RL}.
The consumption of the appliances in a building was scheduled by DQN in \cite{2018-building-RL}.
Then, in \cite{2018-EV-charg-RL}, DQN was used to schedule the EV charging tasks based on the forecasted electricity price information. 
Another category of RL is called policy gradient, that has been applied in \cite{2018-bs-RL-DR,2019-DRL-HVAC}.
Policy gradient was employed for scheduling the consumption of appliances in \cite{2018-bs-RL-DR}.
The similar method was used for managing heating, ventilation, and air conditioning (HVAC) system in the household \cite{2019-DRL-HVAC}; the thermal comfort can be predicted by a Bayesian network while controlling the HVAC system.
However, for both DQN and policy gradient, the action space still needs to be discretized.
Also, with more appliances and households, the scalability of the model-free methods becomes an issue.

In this paper, we propose a distributed learning-based framework for jointly scheduling energy consumption and preserving privacy for the households. 
Specifically, the electricity price is influenced by the total consumption of the households, and therefore the households should work together to reduce the electricity cost. 
The non-cooperative stochastic game is then utilized to model the interactions between the households and the power
grid with the electricity price as an uncertainty factor to all households. 
To search for the NE for the game, a model-free method is adopted, which is different from \cite{2018-Eh-DR-incentive,2018-pil-DR-QOE,2018-DR-dis,2018-DR-uncertainty,2018-DR-dis-shape,2018-DR-dis-EV} that require the price information and consumption ahead of time. 
To further address the discretization issue associated with Q-learning \cite{2015-RL-DR,2019-lu-DR-RL,2018-water-heater-RL-DR,2016-RL-DR-kim,2017-RL-DR,2019-amir-DR-RL-privacy}, DQN \cite{2015-DQN-nature,2018-building-RL,2018-EV-charg-RL}, and policy gradient \cite{2018-bs-RL-DR,2019-DRL-HVAC}, the employed model-free approach is similar to model-based methods that provide a continuous value as the solution without discretizing the action space.
Also, with more appliances and households, the scalability of the model-free methods is an important issue, which is a largely unexplored problem.
To address this issue, a distributed framework is designed to utilize the model-free method, which is motivated by \cite{2014-DR-dis,2018-DR-dis, 2018-DR-dis-shape, 2018-DR-dis-EV}.
When training the parameters of the distributed model-free method, the households may need to exchange their information with each other which is similar as \cite{2014-DR-dis}. 
For privacy concerns, some households are not interested in revealing their consumption profiles to others. 
We therefore introduce a trusted third party is introduced to the proposed framework, such that the privacy of the households can be preserved.

The main contributions of this paper are threefold:
\begin{itemize}

\item The interaction between the power grid and the households is modeled in the form of a stochastic game where the electricity price is the stochastic variable.
Then, we devise a novel model-free approach to search for the NE for the game.

\item We propose a new deep RL approach, which does not require discretization, for scheduling consumption of the appliances in the households.
Also, the privacy of the households is taken into account when utilizing the proposed RL method.

\item We evaluate the performance of the proposed method with a real-world dataset.
The results indicate the proposed model-free method can result in a saving of about $12\%$ in PAR and $11\%$ in load variance.

\end{itemize}

The rest of this paper is organized as follows.
We begin by introducing the system model in Section \ref{sec:system_model_problem_formulation}.
Then, the interactions of the households and the aggregator is formulated as a non-cooperative stochastic game in Section \ref{sec:problem_formulation}.
The way of utilizing reinforcement learning to search for the NE is explained in Section \ref{sec:RL_know}.
Next, utilizing the real-world dataset to evaluate the proposed method as well as the results of the evaluation are provided in Section \ref{sec:simulation}.
Section \ref{sec:conclusion} gives conclusions and suggestions for future work.

\section{System Model}\label{sec:system_model_problem_formulation}

\subsection{System Model}

We consider a power grid with a utility company (i.e., distributed system operator (DSO)), which provides power to a set ${\cal N} = \{ 1, 2, \cdots, N\}$ of $N$ households as shown in Fig. \ref{fig:System_Model}. 
Each household is equipped with an energy consumption controller (ECC), which is responsible to schedule the consumption of the appliances in the household.
ECCs are responsible to schedule the consumption of appliances in a set of equal-length time slots ${\cal T} = \{ 1, 2,\ldots, T \}$.

\begin{figure}
\begin{center}
\resizebox{3.2in}{!}{%
\includegraphics*{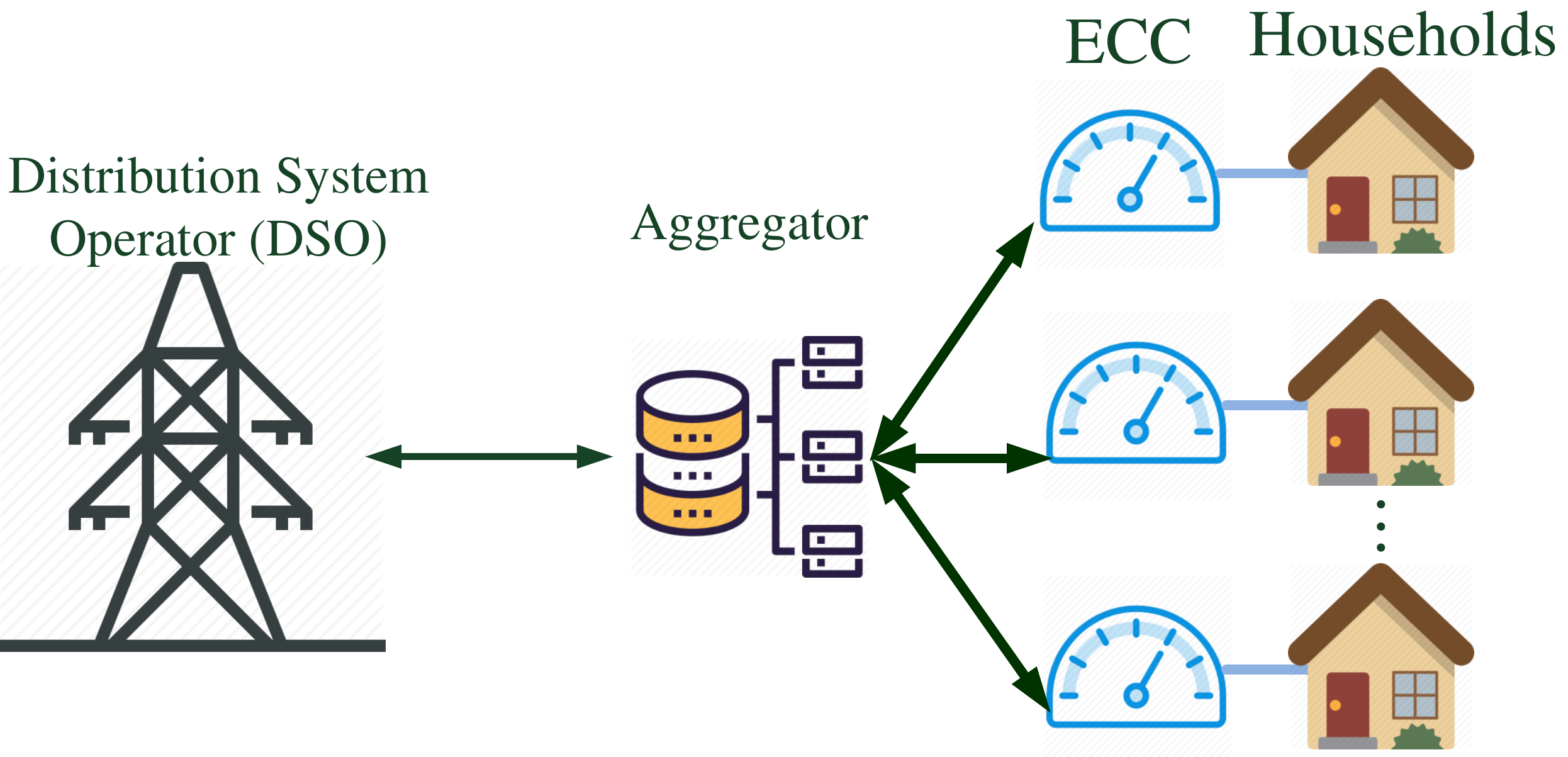} }%
\caption{System model used in this paper.} 
\label{fig:System_Model}
\end{center}
\end{figure}

Let ${\cal A}_{i} = \{ 1, 2, \cdots, A_{i}\}$ denote set of the appliances in household $i$.
In each time slot, an appliance is either awake or asleep.
The appliance is awake when it needs to consume power, and the appliance is regarded as asleep when ECC fulfills its consumption.
If the appliance is awake, it has a state denoted by a tuple as $s_{i, j, t} = \{ t, E_{i, j, t}\}$, where $t$ and $E_{i, j, t}$ represent the current time and required consumption, respectively.
By contrast, the state for the asleep appliance is represented as $s_{i, j, t} = \{ t, 0\}$.

Normally, these appliances can be separated into two categories, which are shiftable and non-shiftable (must-run) appliances.
If an appliance is regarded as the non-shiftable appliance, then the ECC should supply its consumption without any delay when it is awake. 
The ECC can delay the consumption of the shiftable appliances.
Shiftable appliances can further be defined into interruptible and non-interruptible appliances.
For a non-interruptible appliance (e.g., dishwasher or washing machine), the power supply cannot be interrupted once the ECC turns it on.
However, the power supply can be interrupted for an interruptible appliance (e.g., electric vehicle (EV)).
With these categories, we can assign ${\cal A}_{i} = \{ NSA_{i}, NIA_{i}, IA_{i} \}$, where $MSA_{i}, NIA_{i}, IA_{i}$ denote the set of non-shiftable, non-interruptible, and interruptible appliances in the $i$-th household, respectively.
Since the consumption of non-shiftable appliances cannot be modified, the objective of the $i$-th ECC is then to consume an amount of power, $P_{i, t}$, to support the power for appliance $j \in  \{ NIA_{i}, IA_{i} \}$, in every time slot.

For the privacy reason, the households do not reveal their power consumption information to the DSO.
There is an aggregator belonging to a trusted third party to communicate with the DSO for the households.
Therefore, the ECC has limited capability of observing the state of other households in the power grid; that is, the observation of the ECC is limited to its own household.

\subsection{Pricing Model}

The ECC acts as the consumption scheduler for the household.
The objective of the consumption scheduling is to minimize the electricity cost.
 In time slot $t$, the generation cost is used to describe the trend of the electricity price, and therefore the price can be defined as
\begin{equation} \label{eq:gene_cost}
\lambda_{t} = \alpha_{1} L_{t}^{2} + \alpha_{2} L_{t} + \alpha_{3}, 
\end{equation}
where $\alpha_{1}$, $\alpha_{2}$, and $\alpha_{3}$ are the generation coefficients.
Then, $L_{t}$ is the total load in the grid and it is represented as
\begin{equation} \label{eq:net_load}
L_{t} = \sum_{i \in {\cal N}} P_{i, t} +  \sum_{i \in {\cal N}}  \sum_{j \in NSA_{i}} E_{i, j, t}. 
\end{equation}
Besides, the results in \cite{2018-DR-uncertainty} verified that the real-time pricing (RTP) can be implemented in the form of combining the inclining block rate (ICB) mechanism and (\ref{eq:gene_cost}).
Thus, the RTP model used in this paper can be expressed as 
\begin{equation} \label{eq:rtp}
RTP_{t} = \left\{
\begin{array}{ll}
                \lambda_{t},    & 0 \leq  L_{t}  \leq \delta_{1},\\
                \sigma_{1}\lambda_{t} 	& \delta_{1} \leq  L_{t}  \leq \delta_{2},\\
                \sigma_{2}\lambda_{t} 	& L_{t}  > \delta_{2},
\end{array}
        \right.
\end{equation}
where $RTP_{t}$ is the real-time price at time $t$.
$\delta_{1}$ and $\delta_{2}$ indicate the specified net load thresholds.
Then, $\sigma_{1}$ and $\sigma_{2}$ are the threshold pricing constants such that $\sigma_{2} \geq \sigma_{1} \geq 1$.
Since RTP is generated in real-time, the future electricity price is unknown to the ECCs.
The reason of adopting the RTP mechanism is that it can help smooth the load profile.
That is, the electricity price is high if the consumption is high, and therefore the households will be willing to reduce the consumption.
Therefore, the effectiveness of the RTP mechanism can be reflected in the load profile.

\section{Consumption Scheduling Game}\label{sec:problem_formulation}

In order to minimize the electricity cost, each ECC needs to schedule the consumption profile for the household.
The actions of the ECCs are based on the current electricity price information.
However, the electricity price is generated according to the consumption profile of all households.
The electricity price is thus regarded as a stochastic variable for all ECCs.
With this setting, the interaction of the households and the aggregator can be modeled as consumption scheduling game.

\subsection{Game Formulation and Nash Equilibrium}\label{subsec:game_formulation}

When the current electricity price is high, the ECCs shift the consumption to the time slot with a lower electricity price to reduce the electricity cost for the households.
Also, when shifting the consumption, the ECC does not consider how other ECCs react to the high electricity price.
However, the aggregator does not want all ECCs to shift the consumption to the same time slot.
Therefore, the aggregator generates the electricity price again to force ECCs to change their decisions.
That is, the ECCs may face high electricity price again in the future if all ECCs decide to shift the consumption.
At the same time, the current electricity price will decrease.
The ECCs cannot obtain the future electricity price information, and therefore all ECCs should work together to reduce the electricity cost under uncertainty of the electricity prices. 
The interaction among the households and the aggregator can then be captured in the form of a non-cooperative stochastic game.
The main components of the game include:
\begin{itemize}
\item Player: the ECCs in the set of households ${\cal N}$ and the aggregator;
\item State: the states of the appliances in the $i$-th household $\mathbf{s}_{i, t} = \{ s_{i, j, t}| j \in {\cal A}_{i} \}$; 
\item Observation: electricity price and states at time $t$ and $t-1$;
\item Action: power consumption for shiftable appliances in each household $P_{i, t}$; and
\item Payoff function: expected accumulative reward. 
\end{itemize}
At the beginning of time slot $t \in {\cal T}$, the appliances in the households update their states (i.e., whether an appliance is awake or asleep, or it needs more power).
The ECCs submit the consumption to the aggregator, and then the aggregator broadcasts the electricity price according to (\ref{eq:rtp}).
Based on the price information, the ECC of the $i$-th household obtains the observation of the environment denoted by $o_{i, t} = \{ RTP_{t-1}, RTP_{t}, \mathbf{s}_{i, t-1}, \mathbf{s}_{i, t} \} \in {\cal O}_{i}$, where ${\cal O}_{i}$ is the set of possible observations for household $i$.
There exists an action function in each ECC with parameter $\theta_{i}^{\mu}$ for mapping the observation to an action denoted by $P_{i, t} = \pi_{\theta_{i}^{\mu}} (o_{i, t})$.
After deciding $P_{i, t}$, each ECC receives the electricity price again, and then the electricity cost can be represented as 
\begin{equation}\label{eq:r1}
r_{i, t}^{1} = RTP_{t} \times P_{i, t}.
\end{equation}
Other than the electricity cost, the ECC should be aware of a constraint.
That is, the total consumption of the appliances should be the same before and after scheduling.
Another function, $r_{i, t}^{2}$, is then added to address this issue.
The definition of $r_{i, t}^{2}$ is
\begin{equation} \label{eq:r2}
r_{i, t}^{2} = \left\{
\begin{array}{ll}
                0,    & t \neq T  ,\\
                \epsilon_{1} 	& t = T, E_{i, j, t} > 0\\
                \epsilon_{2} 	& t = T, E_{i, j, t} = 0,
\end{array}
        \right.
\end{equation}
where $\epsilon_{2} > 0$ and $\epsilon_{1} <0$.
The reward function is then denoted by $r_{i, t} = r_{i, t}^{2} - r_{i, t}^{1}$.
In each time slot, the ECC aims to minimize the cost; minimizing the cost is the same as maximizing $- r_{i, t}^{1}$.
Then, at the end of the scheduling horizon, the ECCs receive a positive value on $r_{i, t}^{2}$ if all the required consumption of the appliance is fulfilled.
By contrast, $r_{i, t}^{2}$ is negative if the ECCs do not fulfill the requirement.
Thus, the ECCs turn to maximize the reward for the households.
Then, the ECC in the $i$-th household cares about the accumulative reward over a time horizon.  
The accumulative reward of the $i$-th household with discount factor $\gamma_{i}$ at time $k$ can be formulated as 
\begin{equation} \label{eq:expect_cost}
R(o_{i, k}, \pi_{\theta_{i}^{\mu}} (o_{i, k}) ) =  \sum_{t = k}^{T} \gamma_{i}^{t-k} r_{i, t}, 
\end{equation}
where $P_{i, k} = \pi_{\theta_{i}} (o_{i, k}) $.
The discount factor can be used to represent the users' preference.
That is, the users care about the short-term reward if $\gamma$ is close to $0$.
By contrast, the users is foresighted if $\gamma$ is close to $1$. 
The expected accumulative reward can then be extended from (\ref{eq:expect_cost}) as
\begin{align} \label{eq:policy_objective}
J_{i} (\pi_{\theta_{i}^{\mu}}) & = \mathop{\mathbb{E}} \left[ R( o_{i, 1}, P_{i, 1} ) \middle| \pi_{\theta_{i}^{\mu}} \right],  \notag \\
& = \int_{{\cal O}_{i}}  \rho^{\pi_{\theta_{i}}^{\mu}} ( o_{i, 1} )  R \left( o_{i, 1}, \pi_{\theta_{i}^{\mu}} (o_{i, 1})  \right)   d o_{i, t}, \\
& = \mathbb{E}_{ o_{i, t} \sim \rho^{\pi_{\theta_{i}}^{\mu}} } \left[ R \left( o_{i, 1}, \pi_{\theta_{i}^{\mu}} (o_{i, 1}) \right)  \right],  \notag
\end{align}
where $\rho^{\pi_{\theta_{i}^{\mu}}} $ is the distribution of the observation.
In (\ref{eq:policy_objective}), the ECC attempts to maximize the reward function over the time horizon.
For simplicity, we use $P_{i, t}$ instead of $\pi_{\theta_{i}} (o_{i, t})$ later.


After modeling the game, we discuss how each household determines an action function, $\pi_{\theta_{i}} (o_{i, t}) $, in the consumption scheduling game to maximize the expected accumulative reward.
In this context, the Nash equilibrium (NE), which is a solution for the stochastic game, should be searched.

NE is a state of the game where no player can benefit by unilaterally changing strategies, for given strategies of other players.
The NE of the consumption scheduling game, denoted by $\{ \pi_{\theta_{1}^{\mu}}^{*} \cdots \pi_{\theta_{N}^{\mu}}^{*} \}$, which means the ECCs use the optimal parameters for the action function.
Then, for all households, the NE has the following relation:
To obtain the NE, the following optimization problem should be solved
\begin{equation} \label{eq:NE}
J_{i} \left( \pi_{\theta_{i}^{\mu}}^{*}| \pi_{\theta_{-i}^{\mu}}^{*} \right) \geq J_{i} \left(\pi_{\theta_{i}^{\mu}}| \pi_{\theta_{-i}^{\mu}}^{*} \right),
\end{equation}
where $\pi_{\theta_{-i}^{\mu}}^{*}$ is the action functions with optimal parameters in the ECCs except the $i$-th ECC.
The inequality implies that the expected accumulative reward will be lower if the $i$-th ECC does not use optimal parameters for the action function, $\pi_{\theta_{i}^{\mu}}^{*}$.
The proof of the existence of NE is provided in Appendix \ref{subsec:proof_unique_NE}.

To obtain the NE, we need to solve the following optimization problem
\begin{equation} \label{eq:NE_optimization}
\max_{ \theta_{i}^{\mu} } J_{i} \left( \pi_{\theta_{i}^{\mu}}| \pi_{\theta_{-i}^{\mu}} \right)
\end{equation}
In (\ref{eq:NE_optimization}), the $i$-th household needs to know the action functions of other households.
However, the ECCs are not able to observe the states of the appliances in other households.
Therefore, the action functions of other households cannot be obtained.
To address this issue, some model-based methods attempt to model the actions of other households.
However, the performance depends on the accuracy of the models.
Compared to the model-based methods, model-free methods can obtain good performance without modeling the actions of other households.
Hence, we attempt to introduce a model-free method (i.e., policy gradient) to search for the NE.

\section{Proposed Distributed Deep RL Framework} \label{sec:RL_know}

The ECCs need to know the future electricity price to solve (\ref{eq:NE_optimization}).
Also, the actions of other ECCs are required.
Because the ECCs cannot observe the actions of other households and obtain the future electricity price, the expected accumulative reward cannot be obtained directly.
In this context, historical data is utilized.
That is, a model-free method (i.e., policy gradient) is introduced to learn from the past experience.
Moreover, the proposed method should preserve the privacy for the households.
The detailed steps of the proposed method are introduced in the following subsections.

\subsection{Deterministic Policy Gradient} \label{subsec:policy_gradient}

In the consumption scheduling game, there is an actor function mapping the observation to the action with parameter $\theta_{i}^{\mu}$.
In this paper, this function can be represented with a DNN, denoted by actor network.
To solve the optimization problem, the gradient method can be utilized to update the actor network parameter.
However, the differentiation of the objective function, $J_{i} (\pi_{\theta_{i}^{\mu}})$ in (\ref{eq:policy_objective}), cannot be computed directly, because the reward in the future is unknown.
Hence, the gradient cannot be calculated.
Moreover, at time $t$, the actions before the current time are already specified by the ECC, and therefore the accumulative reward is fixed.
Therefore, a one-step bootstrapping function is employed to approximate the expected accumulative reward in the future 
\begin{subequations} \label{eq:action_value_func}
\begin{align}
&  Q^{\theta_{i}^{\mu}} (o_{i, t}, P_{i, t})  = \mathop{\mathbb{E}} \left[ R \left( o_{i, t},  P_{i, t} \right)  \middle| o_{i, t}, P_{i, t}; \pi_{\theta_{i}^{\mu}}\right], \\
& = r_{i, t} +  \gamma \int_{{\cal O}_{i}}  p_{i, t} R \left( o_{i, t+1}, P_{i, t+1} \right)  d o_{i, t+1} , 
\end{align}
\end{subequations}
where $Q^{\theta_{i}^{\mu}} (o_{i, t}, P_{i, t})$ is called action value function and $p_{i, t}$ is the state transition probability from time $t$ to $t+1$.
The action value function states that the value of $R \left( o_{i, t}, P_{i, t}  \right)$ is $Q^{\theta_{i}^{\mu}} (o_{i, t}, P_{i, t})$ given the current state and action.
With the action value function, we can now show how to upgrade $\theta_{i}^{\mu}$.

Suppose $Q^{\theta_{i}^{\mu}} (o_{i, t}, P_{i, t})$ and $\pi_{\theta_{i}^{\mu}} (o_{i, t}) $ are continuous and differentiable.
Then, the policy gradient can be expressed as 
\begin{subequations} \label{eq:DPG}
\begin{align} 
& \nabla_{\theta_{i}^{\mu}}  J_{i} (\pi_{\theta_{i}^{\mu}})  = \nabla_{\theta_{i}^{\mu}} \int_{{\cal O}_{i, t}} Q^{\theta_{i}^{\mu}} (o_{i, t}, P_{i, t}) \rho^{\pi_{\theta_{i}}^{\mu}}  ~ d o_{i, t}, \label{eq:extect_Q}\\
& = \int_{{\cal O}_{i}} \rho^{\pi_{\theta_{i}}^{\mu}}  \nabla_{\theta_{i}^{\mu}} \pi_{\theta_{i}^{\mu}} (o_{i, t})  \nabla_{P_{i, t}}  Q^{\theta_{i}^{\mu}} (o_{i, t}, P_{i, t})   ~ d o_{i, t}, \label{eq:chain_rule}\\
& = \mathbb{E}_{ o_{i, t} \sim \rho^{\pi_{\theta_{i}}^{\mu}} } \left[ \nabla_{\theta_{i}^{\mu}} \pi_{\theta_{i}^{\mu}} (o_{i, t}) \nabla_{P_{i, t}}  Q^{\theta_{i}^{\mu}} (o_{i, t}, P_{i, t}) \right]. \label{eq:DPG_update}
\end{align}
\end{subequations}
The derivation can be obtained in the similar way as policy gradient in \cite{1999-policy-gradient}.
The method in (\ref{eq:DPG}) is denoted by deterministic policy gradient (DPG), which is different from the traditional policy gradient.
The traditional policy gradient outputs a probability distribution function (PDF) over the action space.
However, the output of the consumption scheduling game should only be a scalar.
Thus, DPG \cite{2014-david-dpg} is employed, which does not discretize the action space, in this paper.

In (\ref{eq:DPG_update}), two terms should be determined, $\nabla_{\theta_{i}^{\mu}} \pi_{\theta_{i}^{\mu}} (o_{i, t})$ and $\nabla_{P_{i, t}}  Q^{\theta_{i}^{\mu}} (o_{i, t}, P_{i, t})$.
The former term is easy because the relation between the actor network and the variables is already clear.
However, getting the derivation of the latter term is challenging because there exists an integration in (\ref{eq:action_value_func}).
Thus, another method is required to address this issue.

\subsection{Actor-Critic Method} \label{subsec:actor_critic}

Deriving $Q^{\theta_{i}^{\mu}} (o_{i, t}, P_{i, t})$ is challenging because it contains an integration and the future reward is unknown.
Therefore, a differentiable action value function is substituted in place of $Q^{\theta_{i}^{\mu}} (o_{i, t}, P_{i, t})$.
Normally, there are two alternatives for the substitute in the literature, namely online and offline methods.
The online method is State--action--reward--state--action (SARSA) and the offline method is Q-learning as introduced in \cite{2018-RL-book}.
In this paper, we choose to use the offline method because ECC may record the data locally to train the replaced action value function.
Moreover, similar to DQN, DNN can be utilized to represent the Q-table in Q-learning.
Thus, the action value function can be expressed as 
\begin{equation} \label{eq:Q_network}
 Q^{\theta_{i}^{\mu}} (o_{i, t}, P_{i, t}) \approx  Q^{\theta_{i}^{Q}} (o_{i, t}, P_{i, t}) = \pi_{\theta_{i}^{Q}} (o_{i, t}, P_{i, t}) ,
\end{equation}
where $Q^{\theta_{i}^{Q}} (o_{i, t}, P_{i, t})$ is called a critic, and it is generated by a critic network, $\pi_{\theta_{i}^{Q}} (o_{i, t}, P_{i, t})$, with the parameter, $\theta_{i}^{Q}$.
The objective function of the critic network is 
\begin{equation} \label{eq:Q_obj}
\min_{\theta_{i}^{Q}} \mathop{\mathbb{E}}  \left[ Q^{\theta_{i}^{Q}} (o_{i, t}, P_{i, t}) - y_{i}\right]^{2}
\end{equation} 
$y_{i}$ is the target action value function, and it is generated by
\begin{equation} \label{eq:target_gene}
y_{i} = r_{i, t} + \gamma \hat{Q}^{\theta_{i}^{Q}} (o_{i, t+1}, \hat{P}_{i, t+1}),
\end{equation}
where $\hat{Q}^{\theta_{i}^{Q}} (o_{i, t+1}, \hat{P}_{i, t+1})$ is the target action value for time $t+1$.
Then, $\hat{P}_{i, t+1}$ is obtained by another actor network denoted by target actor network, which is motivated by \cite{2015-DQN-nature}.
The architecture of the actor-critic method is provided in Fig.\ref{fig:actor_critic}.

\begin{figure}
\begin{center}
\resizebox{3in}{!}{%
\includegraphics*{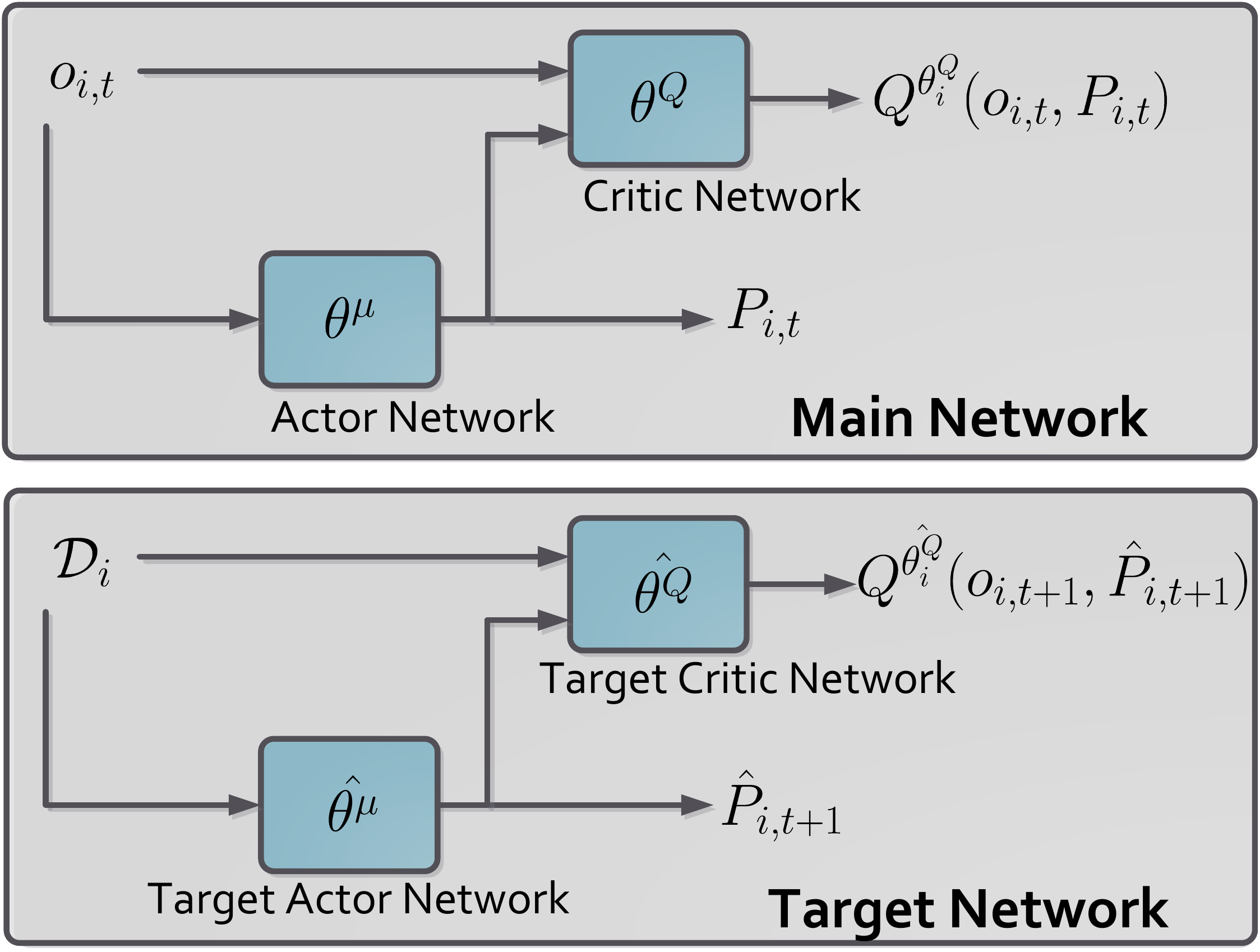} }%
\caption{Actor-critic architecture.} 
\label{fig:actor_critic}
\end{center}
\end{figure}

In order to minimize (\ref{eq:Q_obj}), the observation and the action at the current time slot are not enough because this can introduce bias in the critic network. 
More specifically, the observations at the current time are highly correlated to the observations of previous time slots.
Thus, the training may be unstable.
To overcome this situation, a replay buffer ${\cal D}_{i}$, which stores historical data for the $i$-th ECC, is introduced. 
Then, a subset in ${\cal D}_{i}$ is used to train the critic network, and the number of data in this subset is called batch size, $M$. 
With the critic network, the actor network can be trained.
Also, the critic network can be updated by using the gradient of (\ref{eq:Q_obj}).

\subsection{Centralized Critic Distributed Action}  \label{subsec:cen_critic}

Each ECC can have an actor network and a critic network based on the framework in Fig. \ref{fig:actor_critic}.
However, with a local critic network in the ECC, the parameters of the actor network can be trapped around the local optimum due to the lack of the observations from other households.
To address this problem, we attempt to generate the critic with the observation of all households.
In this context, the aggregator has the global view of the system, and therefore it can generate the critic for each household.
That is, the critic network is moved to the aggregator.
The critic generated by the aggregator can be expressed as 
\begin{equation} \label{eq:Q_at_aggregator}
Q^{\theta_{i}^{Q}} = \pi_{\theta_{i}^{Q}} ( o_{1, t}, \cdots, o_{N, t}, P_{1, t}, \cdots, P_{N, t}  ), 
\end{equation}
where the inputs contain the observations and actions of all households.
In this context, the aggregator is responsible to generate a centralized critic, and then each ECC can still obtain a distributed action.
With this method, the households do not need to exchange $\theta_{i}^{\mu}$ and $P_{i, t}$ with each other.
Thus, the privacy of the households can be preserved.

To implement the centralized critic distributed action, a similar framework as in Fig. \ref{fig:actor_critic} is employed.
That is, two target networks are added, target actor and target critic networks, and they have the same parameters and architecture as actor and critic networks.
The target networks, however, depend on the same parameters we attempt to train.
Hence, this makes the training procedure unstable.
To overcome this situation, the set of parameters is given close to the current parameters but with a time delay.
That is, different from the method in \cite{2015-DQN-nature}, which uses the parameters of current actor and critic networks, the parameters of target networks are updated every $\beta$ ($\beta > 1$) time slots.  
Then, the parameters of target networks can be updated as  
\begin{equation} \label{eq:para_update}
\hat{\theta}_{i}^{\mu} \leftarrow \tau \theta_{i}^{\mu} + (1-\tau) \hat{\theta}_{i}^{\mu},\quad \hat{\theta}_{i}^{Q} \leftarrow \tau \theta_{i}^{Q} + (1-\tau) \hat{\theta}_{i}^{Q},
\end{equation}
where $\hat{\theta}_{i}^{\mu}$ and $\hat{\theta}_{i}^{Q}$ are the parameters of the target actor and target critic networks, respectively.
The relation of all networks is presented in Fig. \ref{fig:ECC_model}.

\begin{figure}
\begin{center}
\resizebox{3.5in}{!}{%
\includegraphics*{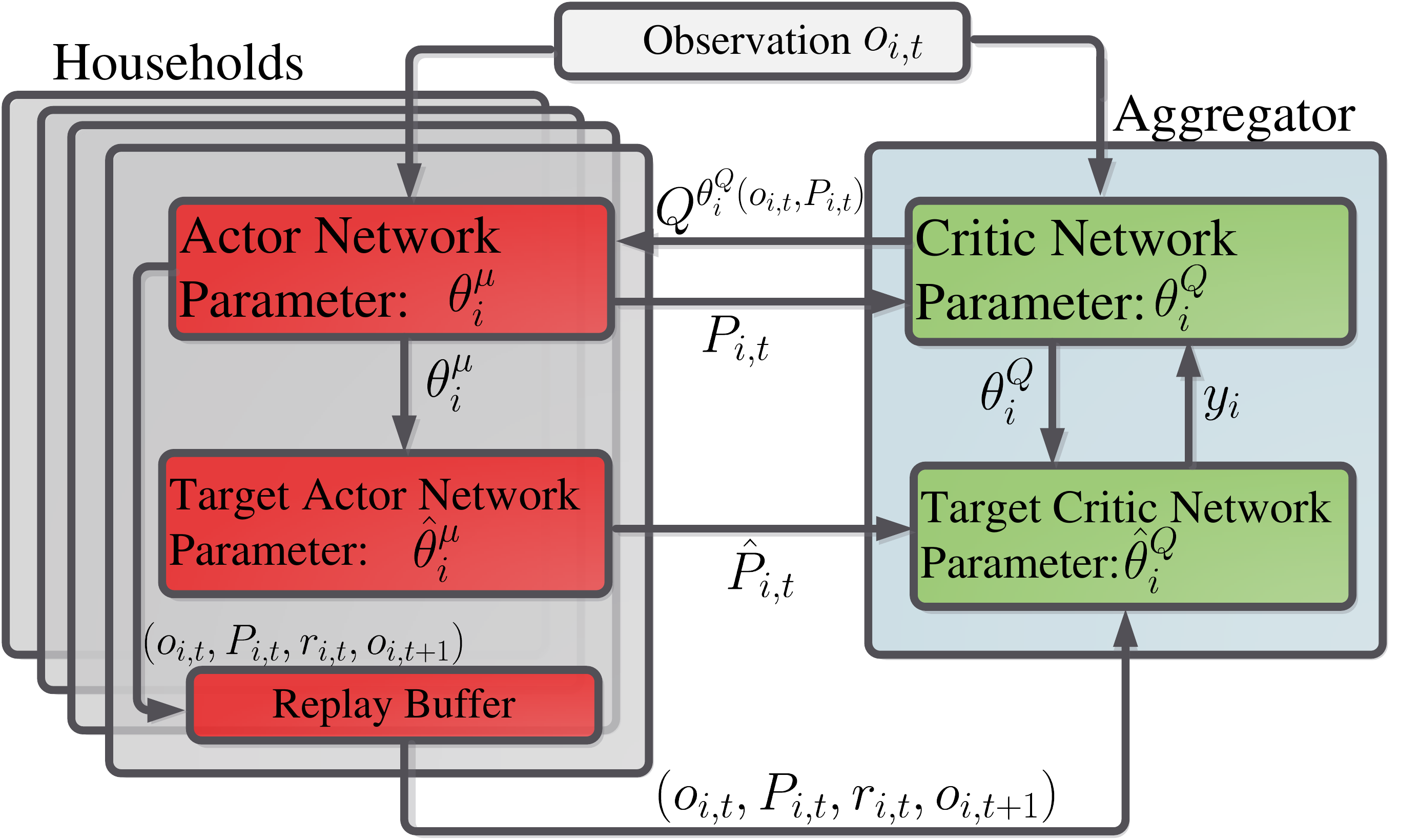} }%
\caption{The architecture of centralized critic distributed action.} 
\label{fig:ECC_model}
\end{center}
\end{figure}

\subsection{Offline Training Online Testing Algorithm}

Based on the techniques introduced from Section \ref{subsec:policy_gradient} to \ref{subsec:cen_critic} and Fig. \ref{fig:ECC_model}, we now summarize the algorithm for training the networks such that the NE of consumption scheduling game can be obtained.
The parameters are trained offline; by doing so, the data in the replay buffer can be used to provide more knowledge to the networks.
In the training phase, the data from $t=1$ to $t=T$ is called an episode.
The training stops when the maximum number of specified episodes is reached.
After the training phase, the trained networks can schedule the consumption of the appliances in real-time.

Before starting the training procedure, the parameters of the networks are initialized randomly.
An additional parameter, $l$, is initialized to $0$. 
Then, the training of the networks is performed by the households interacting with the aggregator.
At the beginning of each time slot, the ECCs first obtain the states of the appliance.
The ECC submits the states to the aggregator, and then receives the electricity price.
By combining the states and the electricity price, the observation can be obtained.
The observation is provided as an input to the actor network, and then $P_{i, t}$ is generated.
During the training, the state space should be explored so that the actor network will not converge to the local optimum.
Thus, a random noise is added to the output of the actor network for exploration as
\begin{equation}
P_{i, t} = \pi_{\theta_{i}^{\mu}} (o_{i, t}) + {\cal N}, 
\end{equation}
where ${\cal N}$ is the exploration noise.
The ECC will decide to draw the power from the grid with the amount of $P_{i, t} + \sum_{j \in NSA_{i}} E_{i, j, t}$.
With the power, the ECC will allocate the power according to the order of non-shiftable, non-interruptible and interruptible appliances.
At the same time, the reward can be obtained through (\ref{eq:r1}) and (\ref{eq:r2}), and then we store $(o_{i, t}, P_{i, t}, r_{i, t}, o_{i, t+1})$ to the replay buffer.
At each time slot, $l$ is increased by $1$.

The parameters of the networks are updated when $l$ reaches $\beta$.
First, $M$ samples from the replay buffer are taken out.
With these samples, $\hat{P}_{i, t+1}$ is generated for the target critic network.
Based on the generated $\hat{P}_{i, t+1}$, the target critic network will generate $y_{i}$ for the critic network.
The critic network can then be updated by minimizing (\ref{eq:Q_obj}).
The critic in the aggregator can broadcast $Q^{\theta_{i}^{Q}}$ to all households based on the updated network.
The action networks of the households can then be updated by using the $M$ samples and (\ref{eq:DPG}).
The parameters of the target actor and the target critic networks are also updated by (\ref{eq:para_update}).
After updating the networks, $l$ is set to $0$ again.
The training procedure operates iteratively until the maximum number of specified episodes, $M_{ep}$, is reached.
Then, the trained networks can be applied to perform real-time consumption scheduling.
The algorithm is denoted by distributed power consumption scheduling (DPCS), which is illustrated in Algorithm \ref{ago:DPCS}.
The convergence of DPCS to NE and the uniqueness of NE can be obtained by starting from the Bellman equation and Bellman error of RL methods as in \cite{2000-RL-proof-ODE,1998-MARL-converge,2002-wang-RL-proof,2009-david-RL-proof}.

\begin{algorithm}
\caption{Distributed Power Consumption Scheduling (DPCS)}
\label{ago:DPCS}
 \DontPrintSemicolon
Initialize actor and critic networks with random weights $\theta_{i}^{\mu}$ and $\theta_{i}^{Q}$ \;
The weight of target actor and target critic networks are assigned as $\hat{\theta}_{i}^{\mu} \leftarrow \theta_{i}^{\mu}$ and $\hat{\theta}_{i}^{Q} \leftarrow \theta_{i}^{Q}$  \;
Initialize the replay buffer ${\cal D}_{i}$ for each ECC and set $l=0$\;
\For{ episode  = $1$ \KwTo $M_{ep}$}{
Obtain the observation $o_{i, t}$ for all households \;
\For{t  = 1 \KwTo $T$ }{	
	Obtain $P_{i, t}$ for all $i$\;
	Execute $P_{i, t}$ and obtain $o_{i, t+1}$ for all $i$\;
	Store $(\mathbf{s}_{i, t}, P_{i, t}, r_{i, t}, \mathbf{s}_{i, t+1})$ to ${\cal D}_{i}$ \;
	Allocate the power to the appliances \;
	$ l = l + 1$ \;
}

\If{$l == \beta$}{
Every ECC samples $M$ data from ${\cal D}_{i}$ to aggregator\;
\For{ECC  = i \KwTo $N$ }{	
	Calculate the $y_{i}$ based on (\ref{eq:target_gene}) and (\ref{eq:Q_at_aggregator})\;
	Update the critic by minimizing \ref{eq:Q_obj}\;
	Update the actor network by (\ref{eq:DPG}); 
}
Update the parameters of target networks (\ref{eq:para_update})\;
$l = 0 $\;
}
}
Use $\!P_{i, t} \!\!=\!\! \pi_{\theta_{i}^{\mu}} \!(\!o_{i, t}\!)$ for real-time consumption scheduling
\end{algorithm}

\subsection{Computational Complexity} \label{subsec:compte_complexity}

The computational complexity of DPCS is now analyzed and compared with the computational complexity of centralized critic centralized action and the computational complexity of the traditional model-based method (e.g., \cite{2019-cui-energy-share}).
Specifically, the authors in \cite{2019-cui-energy-share} design an energy management strategy by utilizing a sample weighted average approximation method to search for the NE of the formulated stochastic game.

Centralized critic centralized action means that there is only one main and one target network; an actor network is responsible for generating the actions for all households.
Since DPCS also uses centralized critic, the difference is in the actor network.
Also, the computational complexity of the model-free methods is equally huge in the training phase.
Therefore, the computational complexity of DPCS is compared with centralized critic centralized action for the real-time control.
The computation of DNN without considering activation functions can be viewed as matrix multiplication.
Then, the computational complexity of DNN can be approximately expressed as ${\cal O} (MKV)$, where $M$, $K$, and $V$ are the number of input nodes, number of hidden layers, and number of output nodes.
The complexity of the centralized action is ${\cal O} (N|o_{i, t}|KN)$, where $|o_{i, t}|$ is the cardinality of $o_{i, t}$; however, it is ${\cal O} (|o_{i, t}|K)$ for DPCS.
When the number of households is huge, DPCS can obtain much lower computational complexity than the centralized critic centralized action.

The computational complexity of \cite{2019-cui-energy-share} is ${\cal O} (n^{3}Nb)$ according to \cite{boyd-cvx-book}, where $n$ is the dimension of the variable and $b$ is the required number of iteration.
That is, the computational complexity of solving the quadratic optimization problem for each household is ${\cal O} (n^{3})$, and then there is a total of $N$ households and needs $b$ rounds of iteration.
On the other hand, the computational complexity of DPCS is ${\cal O} (MM_{ep}|o_{i, t}|K)$ in the training phase.
Thus, the complexity is dominated by the batch size, $M$, and the number of episodes, $M_{ep}$.
Usually, the product of $M$ and $M_{ep}$ is much larger than $n^{3}$.
However, the computational complexity becomes ${\cal O} (|o_{i, t}|K)$ when the ECCs perform the real-time control.
In this case, using the ECCs to perform the real-time control can obtain less computational complexity than \cite{2019-cui-energy-share}.

\section{Numerical Results}\label{sec:simulation}

In this section, the proposed method is evaluated on a large real-world database recorded by Pecan Street Inc. \cite{Load_data}.
Pecan Street Inc. provides a real-world testbed, which includes the data of energy consumption profile with more than $1,000$ households in Texas, USA.
Three scenarios are considered, where one utility company serves $4$, $10$, and $50$ households in the system, respectively. 
There is a total of $10$ types of the appliances described in Table \ref{tb:appliance_setting}, and the household may not have all appliances.
The data from $05/01/2017$ to $10/31/2018$ is used to train the ECCs, and then the consumption profile on $11/05/2018$ is used to evaluate the training results.
The $\alpha_{1}$, $\alpha_{2}$, and $\alpha_{3}$ are set to $0.02$, $0.02$, and $0.5$, respectively.
For the ICB, $\sigma_{1}$ is set to $1.1$ and $\sigma_{2}$ is $1.3$.
Then, $\delta_{1}$ and $\delta_{2}$ are set to $50$ kW and $100$ kW, respectively. 
$\epsilon_{1}$ is set to the product of $-60$ and unfulfilled demand, and $\epsilon_{2}$ is $50$.

Before training the networks, the data from Pecan Street Inc. should be preprocessed.
In the dataset, some households may have more than one specific appliance (e.g., a household may have two or more refrigerators), but some households are not.
We then sum up the consumption profile of these appliances as one appliance.
Also, the data of some households under a certain duration is missing. 
The missing value is synthesized with data of previous days.

The parameters used for training the networks are as follow.
For the DPCS, the actor and critic networks are composed of DNN with $2$ hidden layers, and each layer has $64$ neurons.
The batch size is set to $500$. 
The Rectified Linear Unit (ReLU) is used as the activation function, and the batch normalization is applied in each layer.
Adam is adopted for gradient-based optimization and the learning rate is $0.001$.
The size of the replay buffer ${\cal D}$ is $10^{6}$.
$\gamma_{i}$ is set to $0.99$ and $\tau$ is set to $0.1$.
The maximum episode is set to $1000$ and each episode contains $36$ hours.
Then, $36$ hours are separated to $144$ time slots, which means each time slot is $15$ minutes.
The parameters of the networks are updated every $10$ episodes (i.e., $\beta = 1440$).

To compare the load profile after scheduling the appliances, the PAR metric is utilized, which can be defined as
\begin{equation}
PAR = \log \left( \frac{  L_{k}  }{ \frac{1}{T}    \sum_{t=1}^{T} L_{t} } \right), 
\end{equation}
where $ k = \{ t | \argmax L_{t} \}$.
Then, we compare the electricity cost as
\begin{equation} \label{eq:cost_compare}
\sum_{t=1}^{T} \sum_{i=1}^{N}  RTP_{t} \times P_{i, t}.
\end{equation}
The proposed method is compared with two types of deep deterministic policy gradient (DDPG) \cite{2015-david-ddpg}, which is a popular RL method that can outperform the traditional RL methods (e.g., Q-leaning) and it is also used in \cite{2019-DRL-HVAC}.
The first one is a centralized DDPG framework, denoted by C-DDPG.
That is, the C-DDPG has the all observation of the households and it is responsible for scheduling the consumption of all households.
Another is a distributed framework of DDPG, referred to as D-DDPG.
In D-DDPG, every household has a DDPG network that schedules the consumption of the appliances without considering the actions of other appliances.
Then, C-DDPG and D-DDPG use the same parameters as DPCS.
The method from \cite{2019-cui-energy-share} is the sample weighted average approximation that is referred to as SWAA.
DPCS is also compared with SWAA in the simulation.
The simulations for computation time were conducted with Python $3.7.7$ running on Intel i5-8500B computer with 3.0\,GHz CPU and 16\,GB RAM.

\begin{table} \footnotesize
\begin{center}
\caption{The appliances in the households}\label{tb:appliance_setting}
\begin{tabular}{cc||cc|c|c|c|c|}
\hline
\multicolumn{2}{c||}{Shiftable Appliances}				&  Non-shiftable Appliances\\
Interruptible & 		Non-interruptible		&  \\

\hline
EV  							& Dishwasher      		& Water Heater       \\
Air Conditioner 				& Clothes Washer      	& Lighting       \\
Furnace							& 		      			& Refrigerator       \\
								& 		      			& Kitchen Appliances \\
								& 		      			& Other Appliances      \\
\hline
\end{tabular}
\end{center}
\end{table}

\begin{figure}
\begin{center}
\resizebox{3in}{!}{%
\includegraphics*{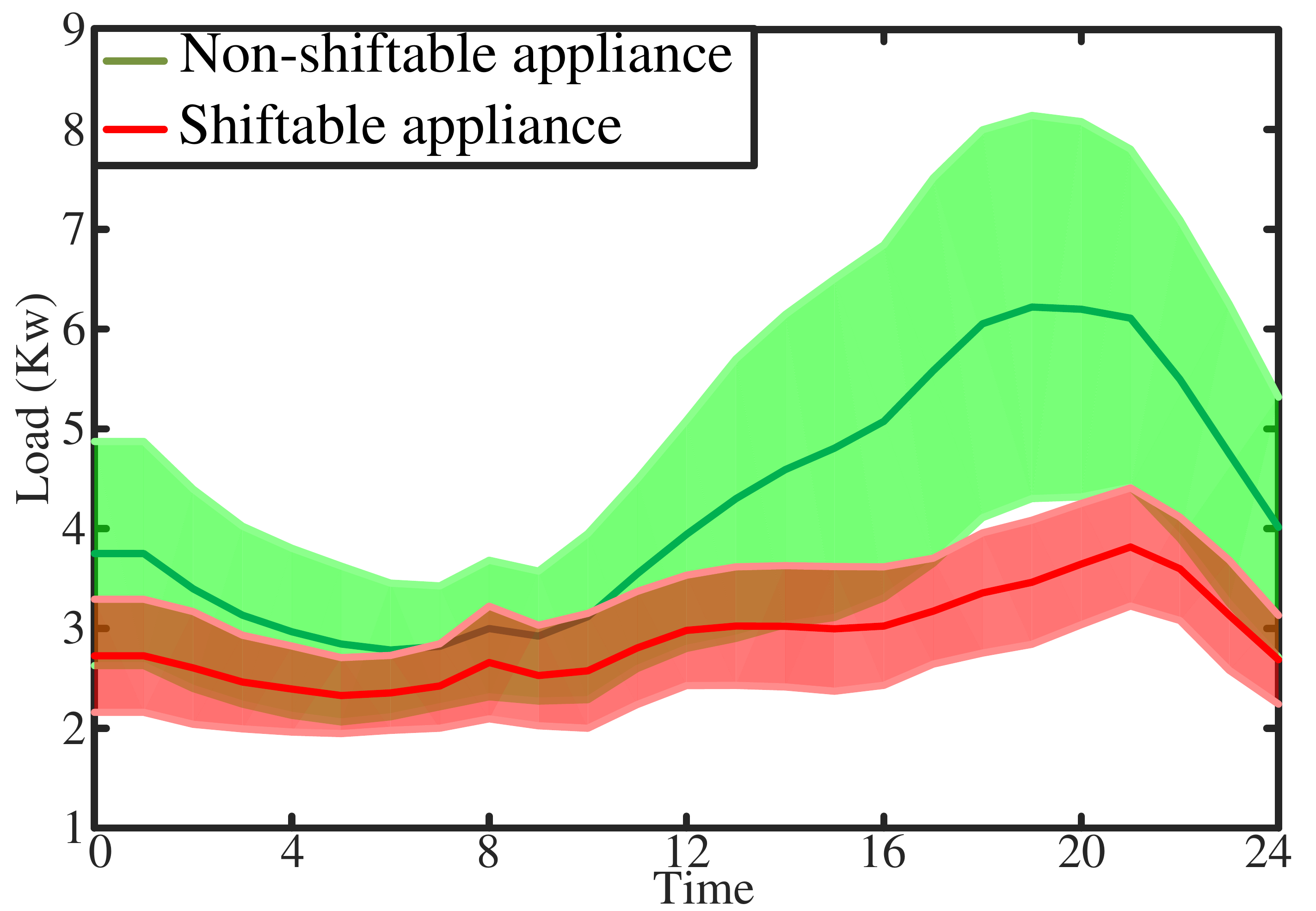} }%
\caption{Load profile of 4 appliances.} 
\label{fig:load_4house_statistic}
\end{center}
\end{figure}

\subsection{Learning Capabilities}

Before assessing the load profile, the training results are provided.
In this part, we only show the results of $4$ households for ease of explaining.
The scenario with $10$ and $50$ households has the same trend and performance.
First, the load curve is shown in Fig. \ref{fig:load_4house_statistic}.
Before $12:00$, the shiftable appliances consume a similar amount as the non-shiftable appliances.
After $12:00$, the consumption of non-shiftable appliance increases drastically, and the load profile of the shiftable appliances remain stable.

After analyzing the load curve, the reward curves versus the number of episodes of three algorithms are shown in Fig. \ref{fig:average_reward}.
The average reward in Fig. \ref{fig:average_reward} represents the trend of the reward during the training phase.
We can observe from Fig. \ref{fig:average_reward} that the average reward is very low at the beginning since the actor networks are initialized with random parameters.
Three algorithms update the policies during the training phase.
For C-DDPG and DPCS, the reward during each episode improves and finally converges to the similar point. 
Then, C-DDPG has all the information of the households and it does not need the help from the aggregator such that it can converge more quickly than DPCS.
In this case, it is proved that the proposed method can converge to the same point as the centralized framework.
However, for the D-DDPG, the ECCs cannot obtain any information from other households such that the ECCs cannot develop good critic network.
Thus, the actor network cannot choose good action, such that the reward keeps decreasing.

\begin{figure}
\begin{center}
\resizebox{3in}{!}{%
\includegraphics*{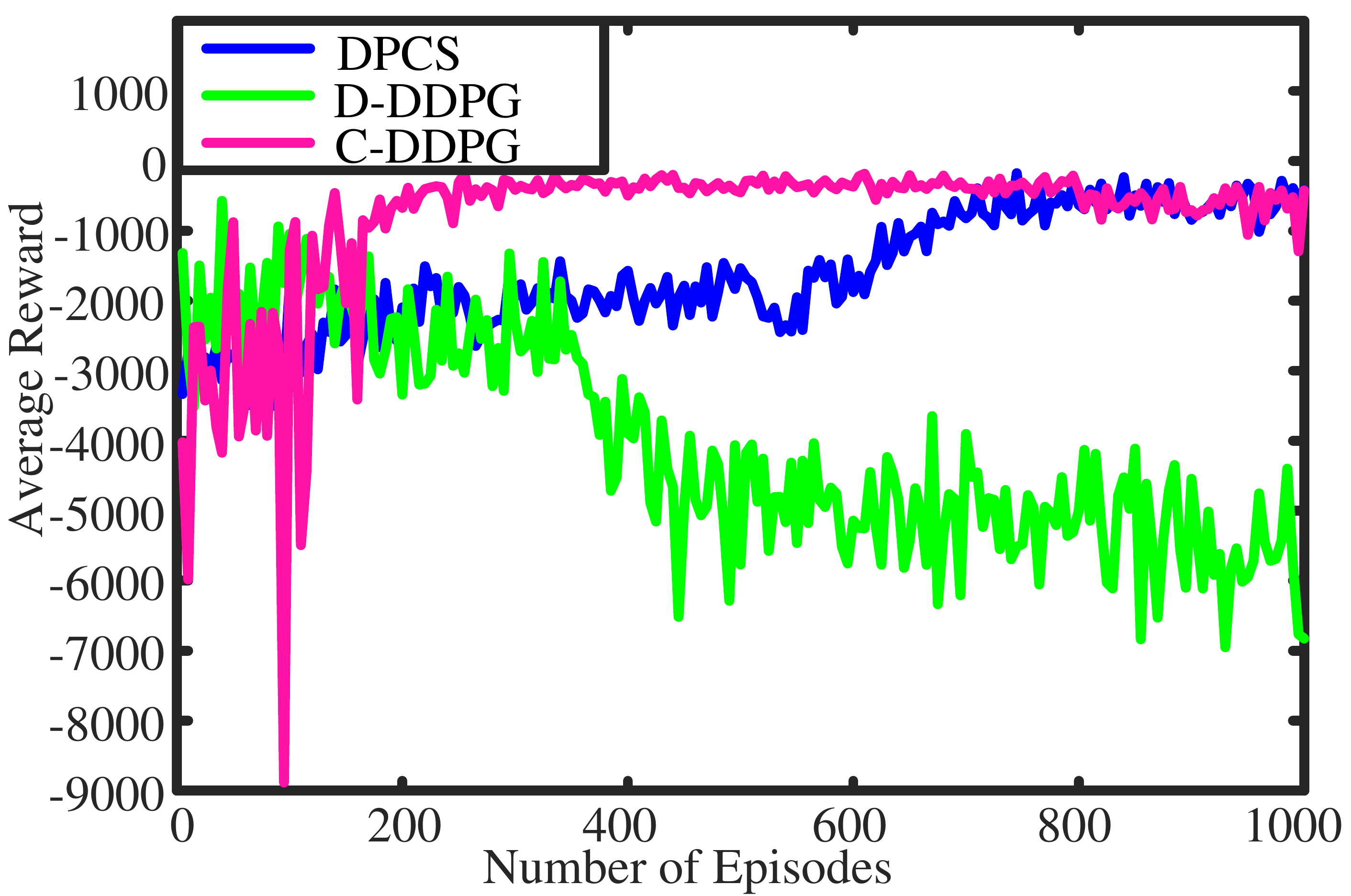} }%
\caption{The average reward for the algorithms.} 
\label{fig:average_reward}
\end{center}
\end{figure}

We now take a closer look at the reward of each household.
The reward curves of each household by applying DPCS are shown in Fig. \ref{fig:individual_reward}.
The objective of DPCS is to maximize the reward function.
For the $3$-rd and the $4$-th household, they obtain poor performance at the beginning of the training; by contrast, household $1$ and $2$ show decent performance from the beginning.
Therefore, during the training phase, households $1$ and $2$ must remain the same performance.
On the other hand, the reward of the $3$-rd and the $4$-th household need to be improved.
According to Fig. \ref{fig:individual_reward}, the proposed method can actually guide the household, which does not have good performance at the beginning, to eventually receive higher reward.

\begin{figure}
\begin{center}
\resizebox{3in}{!}{%
\includegraphics*{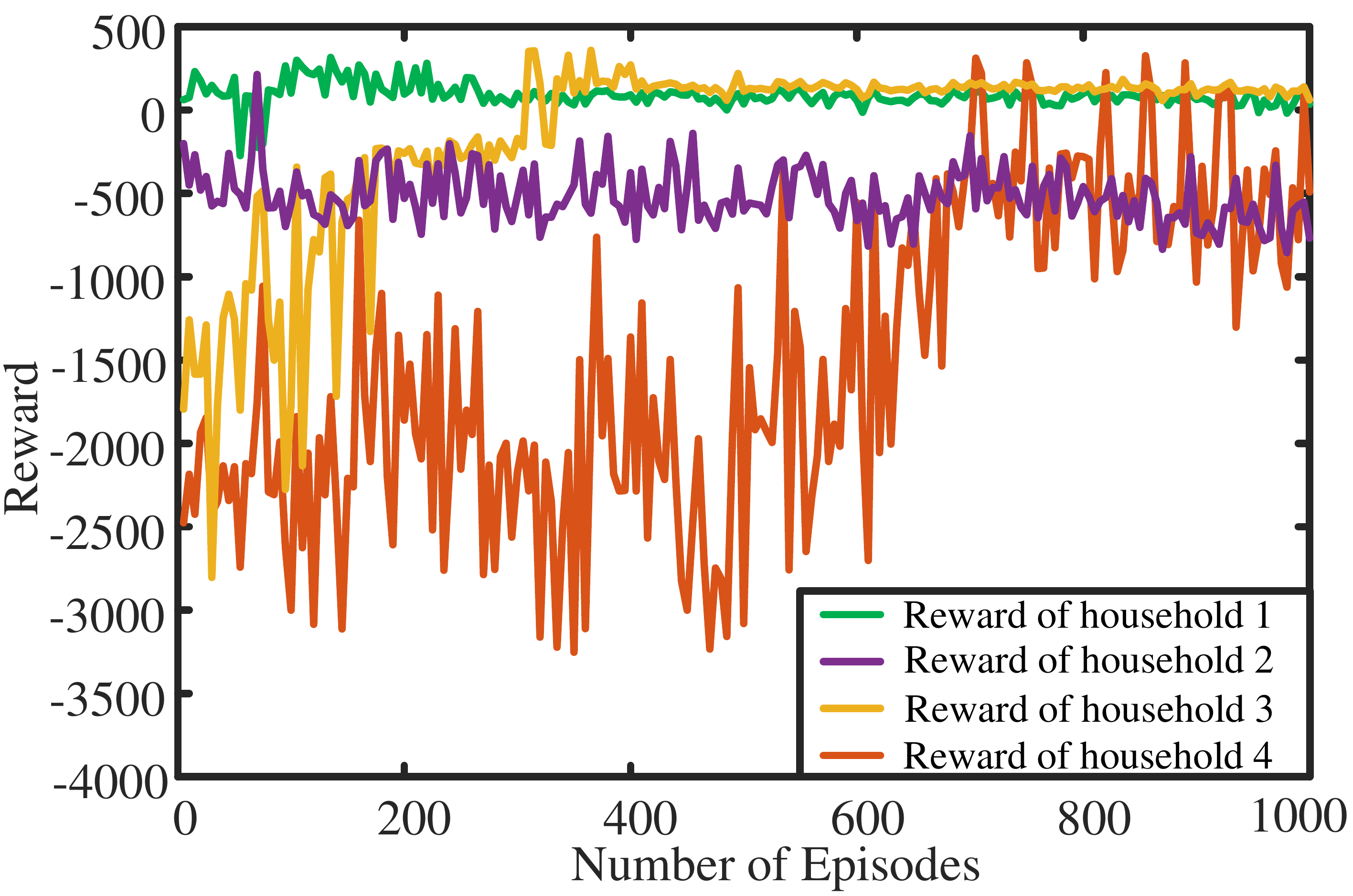} }%
\caption{The individual reward of the households during the training phase.} 
\label{fig:individual_reward}
\end{center}
\end{figure}

Finally, the reward between DPCS and D-DDPG is compared in Table \ref{tb:re_compare}.
At the beginning of the training, the ECCs may have very low reward for both DPCS and D-DDPG.
The ECCs can get higher reward after training with the DPCS for $1000$ episodes.
However, the performance of the ECCs degrades after the training procedure.
This is because the update of action function needs to know the policies of other ECCs.
If the policies of other ECCs cannot be obtained, it can make the actor network trapped in the local optimum or just update randomly.

\begin{table} \small
\begin{center}
\caption{Reward Comparison}
\label{tb:re_compare}
\begin{tabular}{|c|c|c|c|c|c|c|c|c|}
\hline
Method		&  \multicolumn{2}{c|}{DPCS}			& \multicolumn{2}{c|}{D-DDPG} \\
\hline
Episode   & $1$       & $1000$       & $1$       & $1000$   \\
\hline
$1$-st house 	& $-221.67$      & $ 103.55 $      & $ -1969.35 $  & $ -1816.22 $   \\
\hline
$2$-nd house 	& $-1426.35$      & $ -308.32 $      & $ 34.39 $       & $ -2189.60 $   \\
\hline
$3$-rd house 	& $-533.37$       & $ 163.32 $     & $ 147.67 $       & $ -244.11 $   \\
\hline
$4$-th house 	& $-281.43$       & $ 127.65 $     & $ -485.75 $       & $ -275.00 $   \\
\hline
\end{tabular}
\end{center}
\end{table}

\subsection{Load Profile}

To test whether good training results can actually result in promising outcomes for the households, the proposed methods is investigated with the load profile.
We use the trained network to schedule consumption on $11/05/2018$ and then compare the load profile.
The comparing indices are peak load, mean load, variance of the load, PAR, and total electricity cost.
The comparison is provided in Table \ref{tb:load_compare_threealg} and \ref{tb:load_compare}.

\begin{table} \small
\renewcommand{\arraystretch}{1}
\begin{center}
\caption{Load Profile with and without scheduling of DPCS}
\label{tb:load_compare}
\begin{tabular}{|C{2cm}|c|c|c|c|c|c|c|c|}
\hline
	  	&     				  &  Original  		& With scheduling \\
\cline{1-4}
\multirow{5}{*}{4 Households}  & peak (kW) & $ 8.9874 $ & $ 8.3517 $  \\
\cline{2-4}
& mean (kW)      & $ 5.8644 $       & $ 5.7171 $   \\
\cline{2-4}
& Var           & $ 2.4648 $       & $ 2.1378 $   \\
\cline{2-4}
& PAR            &  $ 0.4269 $       & $ 0.3790 $   \\
\cline{2-4}
& Cost (dime)    & $ 205.9080 $       & $ 193.3008 $   \\
\hline 
\hline 
\multirow{5}{*}{10 Households}  & peak (kW) & $19.5222$ &  $18.4471$ \\
\cline{2-4}
& mean (kW)      & $12.9803$       &  $12.9093$  \\
\cline{2-4}
& Var           &  $10.1750$      &  $9.3615$  \\
\cline{2-4}
& PAR            & $0.4081$     &  $0.3570$  \\
\cline{2-4}
& Cost (dime)    & $1483.0473$       & $1369.6979$  \\
\hline 
\hline 
\multirow{5}{*}{50 Households}  & peak (kW) & $97.6112$  &  $92.4692$ \\
\cline{2-4}
& mean (kW)      &  $64.9013$      & $64.6250$   \\
\cline{2-4}
& Var           & $254.3740$        & $226.7057$   \\
\cline{2-4}
& PAR            & $0.4081$   & $0.3597$   \\
\cline{2-4}
& Cost (dime)    & $15813.3927$       &  $15053.4893$ \\
\hline 

\end{tabular}
\end{center}
\end{table}

For the scenario of $4$ households, there is $7.07\%$ reduction on peak load, $11.22\%$ reduction on PAR, and $13.26\%$ on load variance.
For $10$ households, the reduction on peak load, PAR, and load variance are $5.52\%$, $12.52\%$ and the $8.00\%$, respectively.
There is $5.27\%$ reduction on peak load, $12.20\%$ reduction on PAR, and $10.87\%$ reduction on load variance for $50$ households.
According to the simulation results, the proposed method will result in a significant reduction on PAR and lead to a lower variation of the load profile.
With this approach, it will lead to improvement of the stability of the power grid. 
At the same time, the electricity cost for the households can be reduced.

The proposed method can efficiently schedule the power consumption for the households.
However, we do not know the performance of the proposed method compared to C-DDPG even both algorithms can converge to the same point in terms of the average reward.
Also, DPCS should be compared with the model-based approach that is SWAA.
Thus, the comparison of C-DDPG, DPCS, and SWAA is provided in Table \ref{tb:load_compare_threealg}.
In Table \ref{tb:load_compare_threealg}, the computation time is calculated based on the real-time operation in one time slot (15 mins).
According to the results, three algorithms can obtain similar load statistic; it can ensure that DPCS can have the same performance as C-DDPG and SWAA.
However, the load variance in SWAA is less than in DPCS.
This is because SWAA has the information of load profile and the corresponding electricity price for the whole day.
Therefore, it can schedule the load with less variation.
For DPCS, the load profiles of the households undergo quite a bit of fluctuations in the morning as shown in Fig.\ref{fig:load_profile_4house}.
Nonetheless, the performance in terms of the total electricity cost, all three algorithms result in similar values.
On the other hand, DPCS and C-DDPG obtain lower computation time than SWAA which validates our analysis in Section \ref{subsec:compte_complexity}.
Moreover, the proposed method can preserve the privacy for the households; this is not the case with C-DDPG and SWAA.

\begin{table} \small
\renewcommand{\arraystretch}{1.1}
\begin{center}
\caption{Comparison of different scheduling algorithms for $4$ households}
\label{tb:load_compare_threealg}
\begin{tabular}{|c|c|c|c|c|c|c|c|}
\hline
	  	     				  &  C-DDPG  		& DPCS  & SWAA \\
\hline
peak (kW) 		& $ 8.2626 $ & $ 8.3517 $  & $8.2367$ \\
\hline
mean (kW)      & $5.7275 $       & $ 5.7171 $  & $5.7144$ \\
\hline
Var           & $ 2.2880 $       & $ 2.1378 $   & $ 2.0034 $ \\
\hline
PAR            &  $0.3665$       & $ 0.3790 $  & $0.3611$\\
\hline
Cost (dime)    & $ 193.2410 $       & $ 193.3008 $ & $193.1001$ \\
\hline 
\multirow{2}{*}{ \shortstack[l]{Computation time for\\real-time control (s)}}   & \multirow{2}{*}{$1.58$}    & \multirow{2}{*}{$1.53$}   & \multirow{2}{*}{$120.58$}\\
&  & &\\
\hline 
\end{tabular}
\end{center}
\end{table}

In Table \ref{tb:load_compare} and \ref{tb:load_compare_threealg}, one can notice that the reduction of the electricity cost may not be so significant.
This is because of the electricity price and the corresponding load curve.
The electricity price is generated from (\ref{eq:rtp}), which is a quadratic form and highly dependent on the load value. 
Then, non-shiftable appliances consume a huge amount of energy in the evening such that the price in the evening will thus be originally high.
To reduce the cost, the consumption of some shiftable appliances is shifted to the midnight; however, this may increase the electricity price around midnight.
Therefore, the reduction of the electricity cost is mainly on how to shift the consumption of the shiftable appliances in the evening.
Then, in order to virtualize the statistic mentioned above, the load profiles of $4$ households before and after scheduling are presented in Fig. \ref{fig:load_profile_4house}.

We then observe the outcome in Fig. \ref{fig:load_profile_4house}.
In the morning, there may exist peak load for some households (e.g., 8:00 for households $1$ and $2$, or 10:00 for household $3$).
However, this may not cause peak load for the power grid, and therefore the households obtain similar consumption pattern before and after scheduling.
Then, the peak load of the power grid occurs around 20:00 such that the ECCs need to carefully schedule the consumption.
In this context, ECCs $3$ and $4$ decide to decrease the consumption to reduce peak load for the power grid; ECCs $1$ and $2$ remain the same load profile.
This is because households $3$ and $4$ use most of the power for charging EVs, and therefore this portion can be delayed to the midnight.
By doing so, the situation can be avoided where all ECCs decrease the consumption load together when the peak load occurs; the peak load may be shifted to another time slot.

\begin{figure}
\begin{center}
\resizebox{3.55in}{!}{%
\includegraphics*{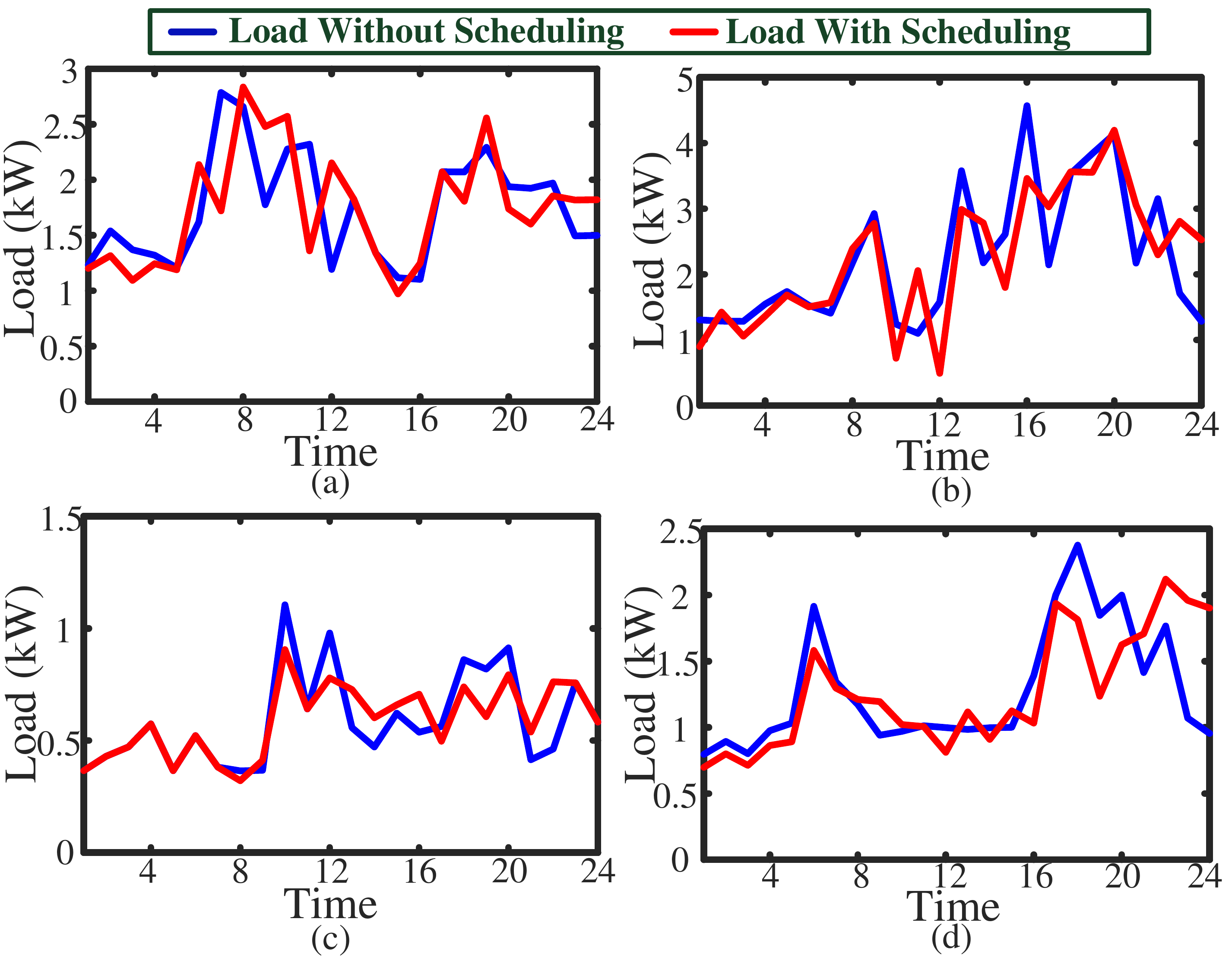} }%
\caption{The load profile of each household. (a) Household 1. (b) Household 2. (c) Household 3. (d) Household 4.} 
\label{fig:load_profile_4house}
\end{center}
\end{figure}

\section{Conclusion}\label{sec:conclusion}

In this paper, a novel model-free framework was proposed for scheduling the consumption profile of the appliances in the households.
The interactions between the households and the aggregator was modeled as a non-cooperative stochastic game.
Then, we applied a model-free method (i.e., policy gradient) to search for the Nash equilibrium (NE).
However, the traditional policy gradient suffers from the discretization issue. 
A method based on deterministic policy gradient (DPG) was then presented to address this issue.
Moreover, the proposed method can preserve the privacy for the households and address the scalability issue.
For verifying the effectiveness of the proposed method, real-world data were employed to evaluate the performance of the proposed method under different number of households in the power grid.
The results revealed that the aggregator introduced in the proposed framework can help the search of the NE for each ECC without the communication between households.
The proposed method was then able to reduce the electricity cost for the households, while also reducing $12\%$ on the PAR for the power grid. 
Moreover, the proposed method can achieve the similar performance as the traditional model-free method and obtain faster real-time control compared to the model-based method.

For future research, we will incorporate the satisfaction level of the households in the proposed method.
That is, the deadline and the comfort level (temperature and humidity) of shifting the consumption should be incorporated.
In addition, it would be interesting to incorporate the feedback of the households to the training procedure. 
In this case, a more customized scheduling policy can be provided for the household.

\appendix

\subsection{Proof of the Existence of NE} \label{subsec:proof_unique_NE}

The following theorem should be proved to ensure the existence of NE of the stochastic game.
\begin{Theorem}\label{the:theorem_exist_NE}
Consider a stochastic game with continuous action and observation spaces. 
If action spaces and observation spaces are compact, payoff functions are continuous, then at least one Nash equilibrium exists.
\end{Theorem}
The definition of compact is provided as follow.
\begin{Definition}\label{def:def_set_compact}
A subset of $\mathbb{R}^{n}$ is compact if the set is closed and bounded.
\end{Definition}
First, the observation space of household $i$, ${\cal O}_{i}$, and the action space of ECC $i$, $A_{i}$, are limited to $[0, \infty)$, and therefore they are compact.
Also, the payoff function is continuous.
Then, we define a function $f_{i}: A_{i} \rightarrow A_{i}$ as $f_{i}(\theta_{i}^{\mu}) = \theta_{i}^{\mu '}$ such that ECC $i$ can obtain higher payoff from $\theta_{i}^{\mu '}$ than from $\theta_{i}^{\mu}$.
Moreover, ECC $i$ can get a benefit of
\begin{equation}
\phi_{i} = \max \left\{ 0, J_{i} \left( \pi_{\theta_{i}^{\mu '}}| \pi_{\theta_{-i}^{\mu}} \right)  - J_{i} \left(  \pi_{\theta_{i}^{\mu}} | \pi_{\theta_{-i}^{\mu}} \right) \right\}, 
\end{equation}
by changing $\theta_{i}^{\mu}$ to $\theta_{i}^{\mu '}$.
Then, Corollary \ref{col:col_fixed_point_sim}, whose proof can be found in \cite{2008-bayesian-NE}, ensures that $f_{i}$ has at least one fixed point.
\begin{Corollary}(Brouwer's fixed point theorem) \label{col:col_fixed_point_sim}
Let $K = \prod \Delta_{m}$, where $\Delta_{m}$ is a simplex, be a simplotope and let $f_{i}: K \rightarrow K$ be continuous. 
Then, $f_{i}$ has a fixed point.
\end{Corollary}
Next, we need to show that the fixed points of $f_{i}$ are Nash equilibria.
If $\theta_{i}^{\mu}$ is a NE, and then $\phi_{i} = 0, \forall i$, making $\theta_{i}^{\mu}$ a fixed point of $f_{i}$.
Conversely, consider $\theta_{i}^{\mu}$ as an arbitrary fixed point of $f_{i}$.
Then, assume that there exists a point, $\theta_{i}^{\mu '}$, that causes $J_{i} \left( \pi_{\theta_{i}^{\mu '}}| \pi_{\theta_{-i}^{\mu}} \right) \geq  J_{i} \left(  \pi_{\theta_{i}^{\mu}} | \pi_{\theta_{-i}^{\mu}} \right)$.
However, $\theta_{i}^{\mu}$ is a fixed point of $f_{i}$; which means $\phi_{i} = 0, \forall i$.
From the definition of $\phi_{i}$, this can only happen if no player can improve its expected payoff by changing $\theta_{i}^{\mu}$.
Thus, $\theta_{i}^{\mu}$ is a NE, and $\theta_{i}^{\mu}$ is equal to $\theta_{i}^{\mu '}$.

{\renewcommand{\baselinestretch}{1}
\begin{footnotesize}
\bibliographystyle{IEEEtran}
\bibliography{References_SmartGrid}
\end{footnotesize}}

\end{document}